# A MONOTONICITY PROPERTY OF RIEMANN'S XI FUNCTION AND A REFORMULATION OF THE RIEMANN HYPOTHESIS

Jonathan Sondow[1] and Cristian Dumitrescu[2]

[1]209 West 97th Street, New York, New York 10025, USA
Email: jsondow@alumni.princeton.edu

[2]119 Young Street, Kitchener, Ontario, N2H4Z3, Canada
Email: cristiand43@gmail.com

### Abstract

We prove that Riemann's xi function is strictly increasing (respectively, strictly decreasing) in modulus along every horizontal half-line in any zero-free, open right (respectively, left) half-plane. A corollary is a reformulation of the Riemann Hypothesis.

## 1. Introduction

The *Riemann zeta function* $\zeta(s)$ is defined as the analytic continuation of the Dirichlet series

$$\zeta(s) = \sum_{n=1}^{\infty} \frac{1}{n^s},$$

which converges if $\Re(s) > 1$. The zeta function is holomorphic in the complex plane, except for a simple pole at $s = 1$. The real zeros of $\zeta(s)$ are $s = -2, -4, -6, \ldots$. Its nonreal zeros lie in the *critical strip* $0 \leq \Re(s) \leq 1$. The *Riemann Hypothesis* asserts that all the nonreal zeros lie on the *critical line* $\Re(s) = 1/2$.

Riemann's *xi function* $\xi(s)$ is defined as the product

$$\xi(s) := \tfrac{1}{2} s(s-1) \pi^{-\tfrac{1}{2}s} \Gamma\!\left(\tfrac{1}{2} s\right) \zeta(s),$$

---

where $\Gamma$ denotes the gamma function. The zero of $s-1$ cancels the pole of $\zeta(s)$, and the real zeros of $s\zeta(s)$ are cancelled by the (simple) poles of $\Gamma(\frac{1}{2}s)$, which never vanishes. Thus, $\xi(s)$ is an entire function whose zeros are the nonreal zeros of $\zeta(s)$ (see [1, p. 80]). The xi function satisfies the remarkable *functional equation*

$$\xi(1-s) = \xi(s).$$

We prove the following monotonicity property of $\xi(s)$. (Throughout this note, *increasing* and *decreasing* will mean strictly so, and a *half-line* will be a half-infinite line not including its endpoint.)

**THEOREM 1.** *The xi function is increasing in modulus along every horizontal half-line lying in any open right half-plane that contains no xi zeros. Similarly, the modulus decreases on each horizontal half-line in any zero-free, open left half-plane.*

For example, since $\xi(s) \neq 0$ outside the critical strip, if $t$ is any fixed number, then $|\xi(\sigma + it)|$ is increasing for $1 < \sigma < \infty$ and decreasing for $-\infty < \sigma < 0$.

In the next section, as a corollary of Theorem 1, we give a reformulation of the Riemann Hypothesis (a slight improvement of [2, Section 13.2, Exercise 1 (e)]). The proof of Theorem 1 is presented in the final section.

## 2. A reformulation of the Riemann Hypothesis

Here is an easy corollary of Theorem 1.

COROLLARY 1. *The following statements are equivalent.*
(i). *If $t$ is any fixed real number, then $|\xi(\sigma + it)|$ is increasing for $1/2 < \sigma < \infty$.*
(ii). *If $t$ is any fixed real number, then $|\xi(\sigma + it)|$ is decreasing for $-\infty < \sigma < 1/2$.*
(iii). *The Riemann Hypothesis is true.*

PROOF. If $|\xi(s)|$ is increasing along a half-line $L$ (or decreasing on $L$), then $\xi(s)$ cannot have a zero on $L$. It follows, using the functional equation, that each of the statements (i) and (ii) implies (iii). Conversely, if (iii) holds, then $\xi(s) \neq 0$ on the right and left open half-planes of the critical line, and Theorem 1 implies (i) and (ii).

## 3. Proof of Theorem 1

We prove the first statement. The second then follows, using the functional equation.

Let $H = H(\sigma_0) = \{s : \Re(s) > \sigma_0\}$ be a zero-free, open right half-plane. Fix a real number $t_0$, and denote by $L = L(\sigma_0, t_0)$ the horizontal half-line

$$L = \{\sigma + it_0 : \sigma > \sigma_0\} \subset H = \{\sigma + it : \sigma > \sigma_0\}.$$

In order to prove that $|\xi(s)|$ is increasing along $L$, we employ the *Hadamard product* representation of the xi function [1, p. 80]:

$$\xi(s) = \tfrac{1}{2} e^{Bs} \prod_\rho \left(1 - \frac{s}{\rho}\right) e^{s/\rho}.$$

Here the product is over all nonreal zeta zeros $\rho$, and $B$ is the negative real number

$$B := \tfrac{1}{2}\log 4\pi - 1 - \tfrac{1}{2}C = -0.023095\ldots,$$

where $C$ is Euler's constant.

We first prove that $|1 - (s/\rho)|$ is increasing on $L$. Since $H = \{s : \Re(s) > \sigma_0\}$ is zero-free and $L \subset H$, we have

$$\Re(\rho) \le \sigma_0 < \Re(s) \quad (s \in L).$$

It follows that the distance $|s - \rho|$ and, hence, the modulus $|1 - (s/\rho)| = |s - \rho||\rho|^{-1}$ are increasing along $L$.

We next show that $|e^{s/\rho}|$ is non-decreasing on $L$. (In fact, $|e^{s/\rho}|$ is increasing on $L$, but we do not need this deeper fact.) Let $\rho = \beta + i\gamma$ denote a nonreal zeta zero. Since $\beta = \Re(\rho) \ge 0$, the modulus

$$\left|e^{s/\rho}\right| = e^{\Re(s/\rho)} = e^{(\beta\sigma + \gamma t_0)/(\beta^2 + \gamma^2)}$$

is non-decreasing along $L$.

It remains to overcome the effect of the Hadamard product factor $e^{Bs}$, which, since $B < 0$, is *decreasing* in modulus on $L$. We use the following alternate interpretation of the constant $B$. First, let $\rho_1, \rho_2, \ldots$ be the zeta zeros with positive imaginary part, and write $\rho_n = \beta_n + i\gamma_n$, for $n \ge 1$. Then $B$ is also given by the formulas [1, p. 82]

$$B = -\sum_{n=1}^{\infty}\left(\frac{1}{\rho_n} + \frac{1}{\overline{\rho}_n}\right) = -2\sum_{n=1}^{\infty}\frac{\beta_n}{\beta_n^2 + \gamma_n^2}.$$

For $N \ge 1$, denote the $N$th partial sum of the series for $-B$ by

$$S_N := \sum_{n=1}^{N} \left( \frac{1}{\rho_n} + \frac{1}{\bar{\rho}_n} \right).$$

Note that $-(B + S_N)$ is positive, and that it approaches zero as $N$ tends to infinity.

Now for $N \geq 2$, let $P_N(s)$ be the finite product

$$P_N(s) := \left(1 - \frac{s}{\bar{\rho}_1}\right) \prod_{n=2}^{N} \left(1 - \frac{s}{\rho_n}\right)\left(1 - \frac{s}{\bar{\rho}_n}\right).$$

Then by combining exponential factors, we can write the Hadamard product as

$$\xi(s) = \tfrac{1}{2} e^{(B+S_N)s} \left(1 - \frac{s}{\rho_1}\right) P_N(s) \prod_{n=N+1}^{\infty} \left(1 - \frac{s}{\rho_n}\right) e^{s/\rho_n} \left(1 - \frac{s}{\bar{\rho}_n}\right) e^{s/\bar{\rho}_n}.$$

From what we have shown about $|1 - (s/\rho)|$ and $|e^{s/\rho}|$, both $P_N(s)$ and the infinite product are increasing in modulus along $L$. To analyze the remaining factors on $L$, set $s = \sigma + it_0$ and define the function

$$f_N(\sigma) := \left| \tfrac{1}{2} e^{(B+S_N)s} \left(1 - \frac{s}{\rho_1}\right) \right|^2 = \tfrac{1}{4} e^{2(B+S_N)\sigma} \frac{(\sigma - \beta_1)^2 + (t_0 - \gamma_1)^2}{\beta_1^2 + \gamma_1^2}.$$

A calculation shows that the derivative $f_N'(\sigma)$ is positive if

$$\frac{\sigma - \beta_1}{(\sigma - \beta_1)^2 + (t_0 - \gamma_1)^2} > -(B + S_N).$$

Now fix $\sigma_1 > \sigma_0$. Since $\sigma_1 - \beta_1 \geq \sigma_1 - \sigma_0 > 0$, and $-(B + S_N) \to 0$ as $N \to \infty$, we can choose $N$ so large that $f_N'(\sigma_1) > 0$. Then $f_N'$ is also positive on some open interval $I$ containing $\sigma_1$. It follows that $f_N(\sigma)$ and, therefore, $|\xi(\sigma + it_0)|$ are increasing for $\sigma \in I$. Since $\sigma_1 (> \sigma_0)$ and $t_0$ are arbitrary, the theorem is proved.